\documentclass[10pt]{article}

\usepackage{dirtytalk}

\usepackage{amsfonts}
\usepackage{amssymb}
\usepackage{amsmath}
\usepackage{enumerate}
\usepackage{amsthm}

\newcommand{\cf}{\mathrm{cf}}

\newcommand{\cov}{\mathrm{cov}}

\newcommand{\pp}{\mathrm{pp}}
\newcommand{\PP}{\mathrm{PP}}
\newcommand{\tcf}{\mathrm{tcf}}

\newtheorem{Th}{\bf THEOREM}[section]

\newtheorem{Pro}[Th]{\bf PROPOSITION}

\newtheorem{fact}[Th]{\bf FACT}

\newtheorem{Cor}[Th]{\bf COROLLARY}

\newtheorem{Obs}[Th]{\bf OBSERVATION}

\theoremstyle{definition} 
 
\theoremstyle{remark}

\theoremstyle{question}

\usepackage[english]{babel}
\makeindex

\title{MEETING, COVERING AND SHELAH'S REVISED GCH}
\author{Pierre MATET}

\date{}

\begin{document}

\maketitle

\renewcommand{\thefootnote}{\arabic{footnote}} 	

\renewcommand{\thefootnote}{}                                
\footnotetext{\textit{Keywords} : Revised GCH Theorem, covering numbers, meeting numbers, Shelah's Strong Hypothesis, club principle}



\vskip 0,7cm

\begin{abstract}  We revisit the application of Shelah's Revised GCH Theorem \cite{SheRGCH} to diamond. We also formulate a generalization of the theorem and prove a small fragment of it. Finally we consider another application of the theorem, to covering numbers of the form $\cov (-, -, -, \omega)$.
\end{abstract}

\bigskip

\section{The Revised GCH Theorem}

\bigskip

Let us start with some definitions.

\medskip

Given a set $S$ and a cardinal $\tau$, we let $P_\tau (S) = \{x \subseteq S : \vert x \vert < \tau \}$.

\medskip

Given four cardinals $\rho_1, \rho_2, \rho_3, \rho_4$ with $\rho_1^+ \geq \rho_2$, $\rho_2^+ \geq \rho_3 \geq \omega$ and $\rho_3 \geq \rho_4 \geq 2$, the {\it covering number} $\cov (\rho_1, \rho_2, \rho_3, \rho_4)$ denotes the least cardinality of any $X \subseteq P_{\rho_2}(\rho_1)$ such that for any $b \in  P_{\rho_3}(\rho_1)$, there is $Q \in  P_{\rho_4}(X)$ with $b \subseteq \bigcup Q$.      

\medskip

Given two infinite cardinals $\chi \leq \lambda$, $\chi$ is a $\lambda$-{\it revision cardinal} if $\cov (\lambda, \chi, \chi, \sigma) \leq \lambda$ for some cardinal $\sigma < \chi$. 

\medskip

A {\it revision cardinal} is an infinite cardinal $\chi$ that is a $\lambda$-revision cardinal for every cardinal $\lambda \geq \chi$.
\medskip

\begin{fact} {\rm (\cite{SheRGCH})}  Every strong limit cardinal $\chi$ is a revision cardinal.
\end{fact}

\medskip

This is actually only one of many (not all equivalent) versions of the \emph{Revised GCH Theorem} (with \say{Revised GCH} often abbreviated as \say{RGCH}) that can be found in \cite{SheRGCH} and other articles of Shelah. In this paper we continue our discussion of the RGCH started in \cite{RGCH}. We would like to cater both for the insiders and the outsiders. Not 	an easy task, considering that some of the former tend to see any mention of cardinal exponentiation as a serious offense. As for the latter, they are aware that for two infinite cardinals $\rho = \cf (\rho) \leq \tau$, $\tau^{< \rho} = \max \{ 2^{< \rho}, \cov (\tau, \rho, \rho, \omega) \}$, so they see the point of studying $\cov (-, -, -, \omega)$, but $\cov (-, -, -, \sigma)$ for $\sigma$ uncountable remains often esoteric to them.

\bigskip

\section{Paradise for all ?}

\bigskip

\subsection{The RGCH for insiders}

\bigskip



To make his message more concrete, Shelah occasionally described the RGCH Theorem as saying that Fact 1.1 holds for $\chi = \gimel_\omega$. He kept looking for \say{sufficient conditions for replacing $\gimel_\omega$ by $\aleph_\omega$} \cite{ShePCF}. A great step forward in this direction was  another version of the RGCH Theorem. Some definitions are in order.

\medskip

Let $X$ be a nonempty set, and $h$ be a function on $X$ with the property that $h (x)$ is a regular cardinal for all $x \in X$. We set $\prod h = \prod_{x \in X} h (x)$. Let $I$ be a (proper) ideal on $X$. For $f, g \in \prod_{x \in X} h (x)$, we let $f <_I g$ if $\{x \in X : f (x) \geq g (x) \} \in I$. 

For a cardinal $\pi$, we set $\tcf (\prod h /I) = \pi$ in case there exists an increasing, cofinal sequence $\vec{f} = \langle f_\alpha : \alpha < \pi \rangle$  in $(\prod h, <_I)$.

\medskip

For a nonempty set $A$ of regular cardinals, we put $\prod A = \prod_{a \in A} a$. For each infinite cardinal $\sigma$, we let ${\rm pcf}_{\sigma{\rm -com}} (A)$ be the collection of all cardinals $\pi$  such that $\pi = \tcf (\prod A /I )$ for some $\sigma$-complete ideal $I$ on $A$. We let ${\rm pcf} (A) = {\rm pcf}_{\omega{\rm -com}} (A)$.

\medskip


Given three infinite cardinals $\tau$, $\rho$ and $\theta$, we let ${\rm Reg} (\tau, \rho, \theta)$ denote the collection of all sets $A$ of regular cardinals such that (1) $0 < \vert A \vert < \tau$, and (2) $\rho < \nu < \theta$ for all $\nu \in A$.

\medskip

Given two infinite cardinals $\chi \leq \lambda$, $\chi$ is $\lambda$-{\it pcf-strong} if for every large enough cardinal $\sigma < \chi$, we have that $\vert \lambda \cap {\rm pcf}_{\sigma{\rm -com}} (A) \vert < \chi$ for all $A \in  {\rm Reg} (\chi, \chi, \lambda)$.

\medskip

\begin{fact} \begin{enumerate}[\rm (i)] 
\item {\rm (\cite{SheRGCH})}  Let $\chi$ be a singular cardinal, and $\lambda$ be a cardinal greater than or equal to $\chi$. Suppose that $\chi$ is $\lambda$-pcf-strong. Then $\chi$ is a $\lambda$-revision cardinal.
\item {\rm (\cite{RGCH})} Let $\chi$ be a limit cardinal of uncountable cofinality, and $\lambda$ be a cardinal greater than or equal to $\chi$. Suppose that there are stationarily many singular cardinals $\nu < \chi$ such that $\nu$ is a $\lambda$-revision cardinal. Then $\chi$ is a $\lambda$-revision cardinal.
\end{enumerate}
\end{fact}

\bigskip

\subsection{The Revision Hypothesis}

\bigskip

We let the \emph{Revision Hypothesis} assert that any limit cardinal is a revision cardinal. By Fact 2.1 (ii), the Revision Hypothesis holds just in case $\chi$ is a $\lambda$-revision cardinal whenever $\chi < \lambda$ are two singular cardinals of cofinality $\omega$.

\medskip

Hence by Fact 2.1 (i), the Revision Hypothesis will hold if $\chi$ is $\lambda$-pcf-strong for any two singular cardinals $\chi < \lambda$ of cofinality $\omega$, which will itself follow if $\vert {\rm pcf} (A) \vert < \vert A \vert^{+ \omega}$ for any set $A$ of regular cardinals with $\min A > \vert A \vert$. 

\medskip

The stronger statement that $\vert {\rm pcf} (A) \vert \leq \vert A \vert$ for any set $A$ of regular cardinals with $\min A > \vert A \vert$ is \emph{Shelah's Medium Hypothesis} (SMH). The failure of SMH (and hence that of the Revision Hypothesis) is known to have large cardinal strength. Gitik \cite{GitikII} has constructed (from large large cardinals) a model with a countable $A$ such that $o.t. ({\rm pcf} (A)) = \omega_1 + 1$, so SMH does not necessarily hold. We do not know whether the Revision Hypothesis may fail.

\medskip

Let us introduce some more notation.

\medskip 

Given three infinite cardinals $\sigma$, $\tau$ and $\theta$ with $\sigma \leq \cf (\theta) < \tau < \theta$, we let $\PP_{\Gamma (\tau, \sigma)} (\theta)$ be the collection of all cardinals $\pi$ such that $\pi = \tcf (\prod h /I )$ for some $h$ and $I$ such that
\begin{itemize}
\item $dom (h)$ is an infinite cardinal less than $\tau$.
\item $h (i)$ is a regular cardinal less than $\theta$ for each $i \in dom (h)$.
\item $\sup \{ h (i) : i \in dom (h) \} = \theta$.
\item $I$ is a $\sigma$-complete ideal on $dom (h)$ such that $\{i \in dom (h) : h (i) < \gamma \} \in I$ for all $\gamma < \theta$.
\end{itemize}

We let $\pp_{\Gamma (\tau, \sigma)} (\theta) = \sup \PP_{\Gamma (\tau, \sigma)} (\theta)$.

\medskip

Given a set $S$ and a cardinal $\rho$, we let $[S]^\rho = \{x \subseteq S : \vert x \vert = \rho \}$.

Given four cardinals $\lambda, \tau, \rho, \sigma$, we define ${\rm equal} (\lambda, \tau, \rho, \sigma)$ (respectively, ${\rm equal} (\lambda, \tau, < \rho, \sigma)$) as follows. If there exists $Z \subseteq P_\tau (\lambda)$ with the property that for any $b$ in $[\lambda]^\rho$ (respectively, $P_\rho (\lambda)$), there is $e \in P_\sigma (Z)$ with $b = \bigcup e$, we let ${\rm equal} (\lambda, \tau, \rho, \sigma) =$ the least size of any such $Z$. Otherwise we let ${\rm equal} (\lambda, \tau, \rho, \sigma) = 2^\lambda$.

\medskip

Given an uncountable limit cardinal $\chi$, we set 

\centerline{$\alpha (\chi, \sigma) = \sup \{\cov (\tau, \nu, \nu, \sigma) : \sigma \leq \nu \leq \tau < \chi\}$}

for each infinite cardinal $\sigma < \chi$.

\medskip

We put $\alpha (\chi) = \min \{ \alpha (\chi, \sigma) : \omega \leq \sigma < \chi\}$.

\medskip

Some consequences of the Revision Hypothesis can be derived from the following. 

\medskip

\begin{fact}  {\rm (\cite{RGCH})} Let $\chi$ be an uncountable limit cardinal, and $\lambda$ be a cardinal greater than or equal to $\chi$. Suppose that $\chi$ is a $\lambda$-revision cardinal. Then the following hold :
\begin{enumerate}[\rm (i)] 
\item Suppose that $\chi$ is a singular cardinal. Then 
\begin{enumerate}[\rm (a)]


\item There is a cardinal $\sigma < \chi$ such that $\pp_{\Gamma (\chi^+, \sigma)} (\theta) \leq \lambda$ for any cardinal $\theta$ with $\sigma \leq \cf (\theta) < \chi < \theta \leq \lambda$.
\item There is a cardinal $\sigma < \chi$ such that ${\rm pcf}_{\sigma{\rm -com}} (A) \subseteq \lambda^+$ for any $A \in  {\rm Reg} (\chi, \chi, \lambda^+)$. 
\item $\cov (\lambda, \chi^+, \chi^+, \sigma) \leq \lambda$ for some cardinal $\sigma < \chi$ (and hence $\chi^+$ is a $\lambda$-revision cardinal).
\item  Suppose that $\cf (\lambda) \not= \cf (\chi)$. Then 
\begin{itemize}
\item There is a cardinal $\sigma < \chi$ such that $\pp_{\Gamma (\chi, \sigma)} (\theta) < \lambda$ for any cardinal $\theta$ with $\sigma \leq \cf (\theta) < \chi < \theta \leq \lambda$.
\item There is a cardinal $\sigma < \chi$ such that ${\rm pcf}_{\sigma{\rm -com}} (A) \subseteq \lambda$ for any $A \in  {\rm Reg} (\chi, \chi, \lambda)$. 
\end{itemize}
\item Suppose that $\alpha (\chi) \leq \lambda$. Then there is a cardinal $\sigma < \chi$ such that for every cardinal $\nu$ with $\sigma \leq \nu < \chi$, $\cov (\lambda, \nu, \nu, \sigma) \leq \lambda$ (and hence $\nu$ is a $\lambda$-revision cardinal).
\item Suppose that  $2^{< \chi} \leq \lambda$. Then
\begin{itemize}
\item ${\rm equal} (\lambda, \chi, \chi, \sigma) \leq \lambda$ for some cardinal $\sigma < \chi$.
\item ${\rm equal} (\lambda, \chi, < \chi, \sigma) \leq \lambda$ for some cardinal $\sigma < \chi$.
\item There is a cardinal $\sigma < \chi$ such that ${\rm equal} (\lambda, \kappa^+, \kappa, \kappa) \leq \lambda$ for every regular cardinal $\kappa$ with $\sigma \leq \kappa < \chi$.\end{itemize}
\end{enumerate}
\item Suppose that $\chi$ is a regular cardinal. Then 
\begin{enumerate}[\rm (1)]


\item There is a cardinal $\sigma < \chi$ such that $\pp_{\Gamma (\chi, \sigma)} (\theta) \leq \lambda$ for any cardinal $\theta$ with $\sigma \leq \cf (\theta) < \chi < \theta \leq \lambda$.
\item There is a cardinal $\sigma < \chi$ such that ${\rm pcf}_{\sigma{\rm -com}} (A) \subseteq \lambda^+$ for any $A \in  {\rm Reg} (\chi, \chi, \lambda^+)$. 
\item  Suppose that $\cf (\lambda) \not= \cf (\chi)$. Then 
\begin{itemize}
\item There is a cardinal $\sigma < \chi$ such that $\pp_{\Gamma (\chi, \sigma)} (\theta) < \lambda$ for any cardinal $\theta$ with $\sigma \leq \cf (\theta) < \chi < \theta \leq \lambda$.
\item There is a cardinal $\sigma < \chi$ such that ${\rm pcf}_{\sigma{\rm -com}} (A) \subseteq \lambda$ for any $A \in  {\rm Reg} (\chi, \chi, \lambda)$. 
\end{itemize}
\item Suppose that $\alpha (\chi) \leq \lambda$. Then there is a cardinal $\sigma < \chi$ such that $\cov (\lambda, \nu, \nu, \sigma) \leq \lambda$ for every cardinal $\nu$ with $\sigma \leq \nu < \chi$.
\item Suppose that  $2^{< \chi} \leq \lambda$. Then
\begin{itemize}
\item ${\rm equal} (\lambda, \chi, < \chi, \sigma) \leq \lambda$ for some cardinal $\sigma < \chi$.
\item There is a cardinal $\sigma < \chi$ such that ${\rm equal} (\lambda, \kappa^+, \kappa, \kappa) \leq \lambda$ for every regular cardinal $\kappa$ with $\sigma \leq \kappa < \chi$.\end{itemize}
\end{enumerate}
\end{enumerate}
\end{fact}

\medskip

Of course denying any one of these consequences will show that the Revision Hypothesis does not always hold. If, on the other hand, the Revision Hypothesis turns out to be true, Fact 2.2 will give an idea of how paradise on earth looks like.

\bigskip

\subsection{Le d\'ebut du d\'ebut}

\bigskip

The Revision Hypothesis does hold in situations when $\lambda$ is (very) close to $\chi$.

\medskip

\begin{fact} {\rm (\cite{RGCH})}
\begin{enumerate}[\rm (i)]
\item Let $\chi$ be an uncountable cardinal. Then $\chi$ is a $\lambda$-revision cardinal for any cardinal $\lambda$ with $\chi \leq \lambda < \chi^{+ \cf (\tau)}$, where $\tau$ equals $\chi$ if $\chi$ is limit cardinal, and the predecessor of $\chi$ otherwise. 
\item Let $\chi$ be a limit cardinal, and $\lambda$ be a cardinal greater than $\chi$. Then $\chi^+$ is a $\lambda$-revision cardinal if and only if $cov (\lambda, \chi^+, \chi^+, \sigma) \leq \lambda$ for some cardinal $\sigma < \chi$.
 \end{enumerate} 
\end{fact}

\begin{Pro}  Let $\chi$ be an uncountable limit cardinal, and $\lambda$ a cardinal with $\chi < \lambda < \chi^{+ \chi}$. Then $\chi$ is $\lambda$-pcf-strong. Furthermore, $\cov (\lambda, \chi^+, \chi^+, \sigma) \leq \lambda$ for some cardinal $\sigma < \chi$.
\end{Pro}

{\bf Proof.} To show that $\chi$ is $\lambda$-pcf-strong, let $A \in {\rm Reg} (\chi, \chi, \lambda)$. Then $A \subseteq B$, and hence ${\rm pcf} (A) \subseteq {\rm pcf} (B)$, where $B$ denotes the collection of all regular cardinals $\mu$ with $\chi < \mu < \lambda$. But $B$ is an interval of regular cardinals, so by a result of Shelah (see Corollary 7.2.7 of \cite{HSW}), $\vert {\rm pcf} (A) \vert \leq \vert {\rm pcf} (B) \vert \leq \vert B \vert^{+ 3} < \chi$.

\medskip

Now assume that $\chi$ is singular. Since $\chi$ is $\lambda$-pcf-strong, it is a $\lambda$-revision cardinal by Theorem 2.1 (i). By Fact 2.2, it follows that $\cov (\lambda, \chi^+, \chi^+, \sigma) \leq \lambda$ for some cardinal $\sigma < \chi$.

\medskip

Finally assume that $\chi$ is regular. By Fact 2.3 (i),  $\chi^+$ is a $\lambda$-revision cardinal (i.e. $\cov (\lambda, \chi^+, \chi^+, \chi) \leq \lambda$). By Fact 2.3 (ii), this implies that $\cov (\lambda, \chi^+, \chi^+, \sigma) \leq \lambda$ for some cardinal $\sigma < \chi$.
\hfill$\square$

\medskip

Let $\chi$ be a singular cardinal, and $\lambda$ a cardinal with $\chi < \lambda < \chi^{+ \chi}$. Then by Fact 2.1 (i) and Proposition 2.4, $\chi$ is a $\lambda$-revision cardinal, and so is $\chi^+$. We will next show that if we make the extra assumption that $\chi$ is not a fixed point of the aleph function, then we also have that $\nu$ is a $\lambda$-revision cardinal for any large enough cardinal $\nu < \chi$. We need some preparation.

\medskip

We start by recalling some properties of covering numbers.

\medskip

\begin{fact}   {\rm (\cite[pp. 85-86]{SheCA}, \cite{LCCN})}  Let $\rho_1, \rho_2, \rho_3$ and $\rho_4$ be four cardinals such that $\rho_1 \geq \rho_2 \geq \rho_3 \geq \omega$ and $\rho_3 \geq \rho_4 \geq 2$. Then the following hold :        
\begin{enumerate}[\rm (i)]
\item  If $\rho_1 = \rho_2$ and either $\cf(\rho_1) < \rho_4$ or $\cf(\rho_1) \geq \rho_3$, then $\cov (\rho_1, \rho_2, \rho_3, \rho_4) = \cf(\rho_1)$.           
\item  If either $\rho_1 > \rho_2$, or $\rho_1 = \rho_2$ and $\rho_4 \leq \cf(\rho_1) < \rho_3$, then $\cov (\rho_1, \rho_2, \rho_3, \rho_4) \geq \rho_1$.        
\item  $\cov (\rho_1, \rho_2, \rho_3, \rho_4) = \cov (\rho_1, \rho_2, \rho_3, \max \{\omega, \rho_4\})$.
\item $\cov (\rho_1^+, \rho_2, \rho_3, \rho_4) = \max \{\rho_1^+, \cov (\rho_1, \rho_2, \rho_3, \rho_4)\}$.
\item  If $\rho_1 > \rho_2$ and $\cf(\rho_1) < \rho_4 = \cf(\rho_4)$, then 

\centerline{$\cov (\rho_1, \rho_2, \rho_3, \rho_4) = \sup \{\cov (\rho, \rho_2, \rho_3, \rho_4) : \rho_2 \leq \rho < \rho_1\}$.}

\item   If $\rho_1$ is a limit cardinal such that $\rho_1 > \rho_2$ and $\cf(\rho_1) \geq \rho_3$, then 

\centerline{$\cov (\rho_1, \rho_2, \rho_3, \rho_4) = \sup \{\cov (\rho, \rho_2, \rho_3, \rho_4) : \rho_2 \leq \rho < \rho_1\}$.}

\item   If $\rho_3 > \rho_4 \geq \omega$, then
 
\centerline{$\cov (\rho_1, \rho_2, \rho_3, \rho_4) = \sup \{\cov (\rho_1, \rho_2, \rho^+, \rho_4) : \rho_4 \leq \rho < \rho_3\}$.}
       
\item   If $\rho_3 \leq \rho_2 = \cf(\rho_2)$, $\omega \leq \rho_4 = \cf(\rho_4)$ and $\rho_1 < \rho_2^{+\rho_4}$, then $\cov (\rho_1, \rho_2, \rho_3, \rho_4) = \rho_1$. 
\item   If $\rho_3 = \cf(\rho_3)$, then either  $\cf(\cov (\rho_1, \rho_2, \rho_3, \rho_4)) < \rho_4$, or  $\cf(\cov (\rho_1, \rho_2, \rho_3, \rho_4)) \geq \rho_3$.
 \end{enumerate}
\end{fact}

\medskip

For two infinite cardinals $\sigma \leq \lambda$, let $u (\sigma, \lambda) = \cov (\lambda, \sigma, \sigma, \omega)$. By Fact 2.5 ((i)-(iii)), 

\centerline{$u (\sigma, \lambda) = \cov (\lambda, \sigma, \sigma, 2) \geq \lambda$.}

\medskip

\begin{fact} {\rm (\cite{Fixed})} Let $\chi$ be an uncountable limit cardinal that is not a fixed point of the aleph function. Then there is an infinite $\eta < \chi$ such that $u (\kappa, \tau) < \chi$ for any regular cardinal $\kappa$ with $\eta < \kappa < \chi$ and any cardinal $\tau$ with $\kappa \leq \tau < \chi$.
\end{fact}

\begin{Obs} Let $\chi$ be an uncountable limit cardinal that is not a fixed point of the aleph function. Then there is an infinite $\eta < \chi$ such that $u (\nu, \tau) \leq \chi$ for any cardinal $\nu$ with $\eta < \nu < \chi$ and any cardinal $\tau$ with $\nu \leq \tau < \chi$.
\end{Obs}

{\bf Proof.} By Fact 2.6, there must be an infinite $\eta < \chi$ such that $u (\kappa, \tau) < \chi$ for any regular cardinal $\kappa$ with $\eta \leq \kappa < \chi$ and any cardinal $\tau$ with $\kappa \leq \tau < \chi$. Now let $\nu$ and $\tau$ be two fixed cardinals such that $\nu$ is singular and $\eta < \nu \leq \tau < \chi$. Then clearly, $P_\nu (\tau) = \bigcup_{\kappa \in X} P_\kappa (\tau)$, where $X$ denotes the set of all regular cardinals $\kappa$ with $\eta < \kappa < \nu$. Hence 

\centerline{$u (\nu, \tau) \leq \sup \{ u (\kappa, \tau) : \kappa \in X \} \leq \chi$.}
\hfill$\square$

\begin{Pro}  Let $\chi$ be an uncountable limit cardinal that is not a fixed point of the aleph function, and $\lambda$ a cardinal with $\chi < \lambda < \chi^{+ \chi}$. Then $\alpha (\chi) \leq \lambda$. Moreover there is a cardinal $\sigma < \chi$ such that for every cardinal $\nu$ with $\sigma \leq \nu \leq \chi^+$, $\cov (\lambda, \nu, \nu, \sigma) \leq \lambda$ (and hence $\nu$ is a $\lambda$-revision cardinal).
\end{Pro}

{\bf Proof.} By Observation 2.7, $\alpha (\chi) \leq \chi$. Now appeal to Proposition 2.4 and Facts 2.1 (i) and 2.2 (i) (e).
\hfill$\square$

\bigskip

\subsection{Shelah's conjectures}

\bigskip

Conjecture 1.11 of \cite{SheWhat} asserts that \say{For every $\mu \geq \aleph_\omega$, for every $\aleph_n < \aleph_\omega$ large enough there is no $\lambda < \mu$ of cofinality $\aleph_n$ such that $\pp_{\Gamma (\aleph_{n+1}, \aleph_n)} (\lambda) > \mu$ (or replace $\aleph_n < \aleph_\omega$ by $\aleph_\alpha < \aleph_{\omega^2}$ or even $\aleph_\alpha < \aleph_{\omega_1}$, or whatever)}. According to \cite{SheDreams} (where 1.11 is misstated), there is \say{a quite reasonable hope} to establish the conjecture. It is pointed out in \cite{SheWhat} that it follows from the conjecture that $\omega_\omega$ is a revision cardinal \say{and this implies $\vert \mathfrak{a} \vert \leq \aleph_0 \implies \vert {\rm pcf} (\mathfrak{a}) \vert \leq \aleph_\omega$, while e.g. $\vert \mathfrak{a} \vert \leq \aleph_{\omega n} \implies \vert {\rm pcf} (\mathfrak{a}) \vert  \leq \aleph_{\omega n + \omega}$ implies} that $\omega_{\omega^2}$ is a revision cardinal.

\medskip


\medskip

In \cite{SheMore} Shelah states two related conjectures \say{which we believe}. One asserts that for any uncountable cardinal $\lambda$, $\kappa^+$ is a $\lambda$-revision cardinal for all but finitely many regular cardinals $\kappa < \lambda$, whereas the other, which is weaker, affirms that for any uncountable strong limit cardinal $\chi$, and any cardinal $\lambda \geq \chi$, $\kappa^+$ is a $\lambda$-revision cardinal for all but finitely many regular cardinals $\kappa < \chi$.


\bigskip

\section{Shelah's Weak Hypothesis}

\bigskip

Let us recall some notation. 

\medskip

For a singular cardinal $\theta$ and a cardinal $\tau$ with $\cf (\theta) < \tau < \theta$, the \emph{pseudopower} $\pp_{< \tau} (\theta)$ is defined as the supremum of the set $X$ of all cardinals $\pi$ for which one may find $A$ and $I$ such that
\begin{itemize}
\item $A$ is a set of regular cardinals smaller than $\theta$ ;
\item $\sup A = \theta$ ;
\item $\vert A \vert < \tau$ ;
\item $I$ is an ideal on $A$ such that $\{A \cap a : a \in A \} \subseteq I$ ;
\item $\pi = \tcf( \prod A /I )$. 
\end{itemize}

\medskip

We let $\pp_\nu (\theta) = \pp_{< \nu^+} (\theta)$ for each cardinal $\nu$ with $\cf (\theta) \leq \nu < \theta$, and $\pp (\theta) = \pp_{\cf (\theta)} (\theta)$.

\medskip

(One version of) \emph{Shelah's Weak Hypothesis} (SWH) asserts that for any cardinal $\lambda$, there are only countably many singular cardinals $\theta$ below $\lambda$ such that $\pp (\theta) \geq \lambda$. SWH implies SMH (\cite{GS01}) and hence the Revision Hypothesis. For more on the SWH see \cite{Git20}.


\bigskip

\section{Meeting numbers}

\bigskip

{\it Shelah's Strong Hypothesis} {\rm (SSH)} asserts that $\pp(\theta) = \theta^+$ for every singular cardinal $\theta$.

\medskip

As the following shows, RGCH and SSH are closely associated.

\medskip

\begin{fact} {\rm (\cite{RGCH})}
The following are equivalent :
 \begin{enumerate}[\rm (i)]
 \item SSH holds.
\item  Let $\chi \leq \lambda$ be two infinite cardinals. Then $\chi$ is a $\lambda$-revision cardinal if and only if $\chi \not= ((\cf (\lambda))^+$.
\item $\omega_1$ is a $\lambda^+$-revision cardinal for every singular cardinal $\lambda$ of cofinality $\omega$.
\end{enumerate}
\end{fact}

\medskip

\medskip

\subsection{Meeting numbers vs. D\v{z}amonja-Shelah numbers}

\medskip

Let us first recall the following notion from \cite{DS}.

\medskip

Let $\rho \leq \lambda$ be two infinite cardinals. For $f : \lambda \rightarrow \rho$, let $ds (f)$ denote the least size of any $Z \subseteq [\lambda]^\rho$ with the property that whenever $B \subseteq \lambda$ is such that $\vert B \cap f^{- 1} (\{\alpha\}) \vert = \lambda$ for all $\alpha < \rho$, there is $z \in Z$ such that $\vert f``(z \cap B) \vert = \rho$.

We define the {\it D\v{z}amonja-Shelah number} $DS (\rho, \lambda)$ by : 

\centerline{$DS (\rho, \lambda) = \sup \{ ds (f) : f : \lambda \rightarrow \rho\}$.}

\medskip

Given two infinite cardinals $\rho \leq \lambda$, the {\it meeting number} $m (\rho, \lambda)$ (respectively, the {\it density number} $d (\rho, \lambda)$) 
denotes the least cardinality of any $Q \subseteq [\lambda]^\rho$ with the property that for any $b \in [\lambda]^\rho$, there is $q \in Q$ with $\vert b \cap q \vert = \rho$ (respectively, $q \subseteq b$).

\medskip

Obviously, $m (\rho, \lambda) \leq d (\rho, \lambda)$. It is known \cite{SSH} that  
\begin{itemize}
\item $\lambda \leq m (\rho, \lambda)$ in case $\rho < \lambda$.
\item $m (\rho, \lambda^+) = \max \{ \lambda^+, m (\rho, \lambda)\}$.
\end{itemize}

\medskip

\begin{fact}  \begin{enumerate}[\rm (i)]
\item {\rm (\cite{GR})} The following are equivalent :
 \begin{enumerate}[\rm (1)]
 \item SSH holds.
\item $DS (\cf (\tau), \tau^+) \leq \tau^+$ for every singular cardinal $\tau$.
\end{enumerate}
\item {\rm (\cite{SSH})} The following are equivalent :
 \begin{enumerate}[\rm (1)]
 \item SSH holds.
\item  Given two infinite cardinals $\rho \leq \lambda$, $m (\rho, \lambda)$ equals $\lambda$ if $\cf (\lambda) \not= \cf (\rho)$, and $\lambda^+$ otherwise. 
\item $m (\omega, \lambda) = \lambda^+$ for every singular cardinal $\lambda$ of cofinality $\omega$.
\end{enumerate}
\end{enumerate}
\end{fact}

\medskip

It turns out that Fact 4.2 does not give two reformulations of SSH, but just one. The following has been independently obtained by Ziemek Kostana and Assaf Rinot.

\medskip

\begin{Obs} Let $\rho \leq \lambda$ be two infinite cardinals. Then $DS (\rho, \lambda) = m (\rho, \lambda)$.
\end{Obs}

{\bf Proof.} $\leq$ : Select $Q \subseteq [\lambda]^\rho$ with $\vert Q \vert = m (\rho, \lambda)$ such that for any $b \in [\lambda]^\rho$, there is $q \in Q$ with $\vert b \cap q \vert = \rho$. Given $f : \lambda \rightarrow \rho$, and $B \subseteq \lambda$ such that $\vert B \cap f^{- 1} (\alpha) \vert = \lambda$ for all $\alpha < \rho$, pick $g \in \prod_{\alpha < \rho} (B \cap f^{- 1} (\alpha))$. There must be $q \in Q$ such that $\vert ran (g) \cap q \vert = \rho$. Then clearly, $\vert f``(q \cap B) \vert = \rho$. Thus, $ds (f) \leq m (\rho, \lambda)$.

\bigskip

$\geq$ : Pick a bijection $k : \lambda \times \lambda \times \rho \rightarrow \lambda$, and define $f : \lambda \rightarrow \rho$ by letting $f^{- 1} (\{\alpha\}) = \{ k (\alpha, \beta, \xi) : \beta, \xi < \lambda\}$ for all $\alpha < \rho$. Select $Z \subseteq [\lambda]^\rho$ such that $\vert Z \vert \leq DS (\rho, \lambda)$, and with the property that whenever $B \subseteq \lambda$ is such that $\vert B \cap f^{- 1} (\alpha) \vert = \lambda$ for all $\alpha < \rho$, there is $z \in Z$ such that $\vert f``(z \cap B) \vert = \rho$. For $z \in Z$, put 

\centerline{$q_z = \bigcup_{\alpha < \rho} \{ \xi < \lambda : \exists \beta < \lambda (k (\alpha, \beta, \xi) \in z)\}$.}

Notice that $\vert q_z \vert \leq \rho$. Given $b \in [\lambda]^\rho$, let $\langle \xi_\eta : \eta < \rho \rangle$ be a one-to-one enumeration of $b$, and set $B = \bigcup_{\alpha < \rho} \{k (\alpha, \beta, \xi_\alpha) : \beta < \lambda\}$. There must be $z \in Z$ such that $\vert f``(z \cap B) \vert = \rho$. Then $\vert b \cap q_z \vert = \rho$. Thus, $ds (f) \geq m (\rho, \lambda)$.
\hfill$\square$  



\medskip

\subsection{Meeting numbers and the RGCH}

\medskip

In view of Facts 4.1 and 4.2 (ii), one can expect consequences of the RGCH in terms of meeting numbers.

\medskip

\begin{Obs} Let $\sigma$, $\rho$, $\lambda$ be three infinite cardinals with $\sigma \leq \cf (\rho)$. Then $m (\rho, \lambda) \leq \cov (\lambda, \rho^+, \rho^+, \sigma)$.
\end{Obs}

{\bf Proof.}  Select $Z \subseteq P_{\rho^+}(\lambda)$ with $\vert Z \vert = \cov (\lambda, \rho^+, \rho^+, \sigma)$ such that for any $b \in [\lambda]^\rho$, there is $e_b \in P_\sigma (Z)$ with $b \subseteq \bigcup e_b$. Then clearly, $\vert b \cap z \vert = \rho$ for some $z \in e_b$.
\hfill$\square$  

\begin{Cor} Let $\mu$ be a regular cardinal, and $\lambda$ be a cardinal greater than $\mu$ such that $\mu^+$ is a $\lambda$-revision cardinal. Then $m (\mu, \lambda) \leq \lambda$.
\end{Cor}


\begin{Cor} Let $\chi$ be an uncountable limit cardinal, and $\lambda$ be a cardinal greater than or equal to $\chi$. Suppose that there is a cardinal $\sigma < \chi$ such that $\cov (\lambda, \nu, \nu, \sigma) \leq \lambda$ for every cardinal $\nu$ with $\sigma \leq \nu < \chi$. Then $m (\rho, \lambda) \leq \lambda$ for every cardinal $\rho < \chi$ with $\sigma \leq \cf (\rho)$.
\end{Cor}

\medskip

\subsection{Sure bet}

\medskip

In \cite{RGCH} and the present paper we have seen several equivalent formulations of SSH that each corresponds to a version of the RGCH, so it is more than likely that the topological form of SSH considered by Rinot in \cite{RinotSSH} (where SSH is shown to be equivalent to the statement that $\vert X \vert = \kappa$ whenever $(X, \tau)$ is a first-countable space of density the regular cardinal $\kappa$ with the property that any separable subspace of $X$ has size at most $\kappa$) can be handled in the same way. We leave it to the reader to formulate and prove the corresponding theorem.

\medskip

\subsection{Clubs are a (modern) girl's best friend}

\medskip

One of the main applications considered in \cite{SheRGCH} is to diamond. Let us first recall some notation.

\medskip

Let $\tau$ be a regular uncountable cardinal. Given a $\tau$-complete ideal $J$ on $\tau$, the {\it diamond principle} $\diamondsuit_\tau [J]$ asserts the existence of $s_\alpha \subseteq \alpha$ for $\alpha < \tau$ such that $\{\alpha < \tau : s_\alpha = A \cap \alpha\}$ lies in $J^+$ for every $A \subseteq \tau$. The {\it diamond star principle} $\diamondsuit_\tau^\ast [J]$ asserts the existence of $t^{i}_\alpha \subseteq \alpha$ for $i < \alpha < \tau$ such that $\{\alpha < \tau : \exists i < \alpha (t^{i}_\alpha = A \cap \alpha)\}$ lies in $J^\ast$ for every $A \subseteq \tau$.

\medskip

We let $NS_\tau$ denote the nonstationary ideal on $\tau$.

\medskip

$\diamondsuit_\tau [NS_\tau]$ is abbreviated as $\diamondsuit_\tau$.

\medskip

For a regular cardinal $\mu$, $E_\mu^\tau$ denotes the set of all $\alpha \in \tau$ with $\cf (\alpha) = \mu$.

\medskip

In the introduction of \cite{SheRGCH}, the application is described as follows : \say{we show that for $\lambda \geq \gimel_\omega$, $2^\lambda = \lambda^+$ is equivalent to $\diamondsuit_{\lambda^+}$} (this is Claim 3.2 (a) in \cite{SheRGCH}). In this formulation, the result has been superseded : by a later result of Shelah \cite{She10}, for every cardinal $\lambda \geq \omega_1$ with $2^\lambda = \lambda^+$, $\diamondsuit_{\lambda^+} [NS_{\lambda^+}]$ holds, and in fact so does $\diamondsuit_{\lambda^+} [J]$ holds whenever $J$ is a $\lambda^+$-complete ideal on $\lambda^+$ extending $NS_{\lambda^+} \vert E^{\lambda^+}_\mu$ for a regular cardinal $\mu$ less than $\lambda^+$ with $\mu \not= \cf (\lambda)$. However immediately after the proof of Claim 3.2 there is the remark that \say{Note that we actually proved also}, followed by the statement of Claim 3.6, where (1) gives the following :

\medskip





\begin{fact} {\rm (\cite{SheRGCH})} Let $\chi$ be an uncountable strong limit cardinal, and $\lambda$ a cardinal greater than or equal to $\chi$ such that $2^\lambda = \lambda^+$. Then there is a cardinal $\sigma < \chi$ such that $\diamondsuit_{\lambda^+}^\ast [NS_{\lambda^+} \vert E^{\lambda^+}_\mu]$ holds for any regular cardinal $\mu$ with $\sigma \leq \mu < \chi$.
\end{fact}


\medskip

Now diamond is closely associated with cardinal exponentiation, so the faithful tend to see it as an outdated relic of the bygone pre-pcf theory days. In their opinion priority should be given to club which they consider to be the heart and soul of diamond. Let us recall some versions of this other principle.

\medskip

Let $\tau$ be a regular uncountable cardinal $\tau$, and $J$ be a $\tau$-complete ideal on $\tau$. For a cardinal $\eta < \tau$, $\clubsuit^{cof \slash \eta, \ast}_\tau [J]$ asserts the existence of $B^{i}_\delta \in P_\eta (\delta)$ for $i < \delta < \tau$ such that for any $A \in [\tau]^\tau$,

\centerline{$\{ \delta < \tau : \exists i < \delta (\sup (A \cap B^{i}_\delta) = \delta)\} \in J^\ast$.}



\medskip

$\clubsuit_\tau^\ast [J]$ asserts the existence of $s^{i}_\alpha \subseteq \alpha$ with $\sup s^{i}_\alpha = \alpha$ for each infinite limit ordinal $\alpha < \tau$ and each $i < \alpha$  such that $\{ \alpha < \tau : \exists i < \alpha (s^{i}_\alpha \subseteq A)\} \in J^\ast$ for all $A \in [\tau]^\tau$.

\medskip

\begin{fact} \begin{enumerate}[\rm (i)] 
\item  {\rm (\cite{SSH})} Suppose that $\tau = \nu^+$, and let $\mu$ be a regular cardinal less than $\nu$ such that $m (\mu, \nu) = \nu$. Then $\clubsuit^{cof \slash {\mu^+}, \ast}_\tau [NS_\tau \vert E^\tau_\mu]$ holds.
\item {\rm (\cite{Towers})} Let $\mu$, $\tau$ be two regular cardinals with $d (\mu, \mu) < \tau$. Suppose that $\clubsuit^{cof \slash {\mu^+}, \ast}_\tau [NS_\tau \vert E^\tau_\mu]$ holds. Then so does $\clubsuit^\ast_\tau [NS_\tau \vert E^\tau_\mu]$.
\item {\rm (\cite{Towers})} Let $\mu < \tau$ be two regular cardinals. Then $\diamondsuit^\ast_\tau [NS_\tau \vert E^\tau_\mu]$ holds if and only if $2^{< \tau} = \tau$ and $\clubsuit^\ast_\tau [NS_\tau \vert E^\tau_\mu]$ holds.
\end{enumerate}
\end{fact}

\begin{Pro} Let $\mu, \lambda$ be two infinite cardinals with $\cf (\mu) = \mu < \lambda$. Suppose that $\mu^+$ is a $\lambda$-revision cardinal. Then 
\begin{itemize}
\item $\clubsuit^{cof \slash {\mu^+}, \ast}_{\lambda^+} [NS_{\lambda^+} \vert E^{\lambda^+}_\mu]$ holds.
\item If $d (\mu, \mu) \leq \lambda$, then $\clubsuit^\ast_{\lambda^+} [NS_{\lambda^+} \vert E^{\lambda^+}_\mu]$ holds. 
\item If $d (\mu, \mu) \leq \lambda$ and $2^\lambda = \lambda^+$, then $\diamondsuit^\ast_{\lambda^+} [NS_{\lambda^+} \vert E^{\lambda^+}_\mu]$ holds.
\end{itemize}
\end{Pro}

{\bf Proof.} By Corollary 4.5 and Fact 4.8 (i).
\hfill$\square$




\begin{Cor} 
\begin{enumerate}[\rm (i)] 
\item Assume SSH, and let $\mu, \lambda$ be two infinite cardinals with $\cf (\mu) = \mu < \lambda$ and $\mu \not= \cf (\lambda)$. Then 
\begin{itemize}
\item $\clubsuit^{cof \slash {\mu^+}, \ast}_{\lambda^+} [NS_{\lambda^+} \vert E^{\lambda^+}_\mu]$ holds.
\item If $d (\mu, \mu) \leq \lambda$, then $\clubsuit^\ast_{\lambda^+} [NS_{\lambda^+} \vert E^{\lambda^+}_\mu]$ holds. 
\item If $d (\mu, \mu) \leq \lambda$ and $2^\lambda = \lambda^+$, then $\diamondsuit^\ast_{\lambda^+} [NS_{\lambda^+} \vert E^{\lambda^+}_\mu]$ holds.
\end{itemize}
\item Let $\chi$ be an uncountable strong limit cardinal, and $\lambda$ a cardinal greater than or equal to $\chi$. Then there is a cardinal $\sigma < \chi$ such that $\clubsuit_{\lambda^+}^\ast [NS_{\lambda^+} \vert E^{\lambda^+}_\mu]$ holds for any regular cardinal $\mu$ with $\sigma \leq \mu < \chi$.
\end{enumerate}
\end{Cor}

{\bf Proof.} (i) : Use Fact 4.1.

\medskip

(ii) : Use Facts 1.1, 2.2 ((i) (e) and (ii) (e)).
\hfill$\square$

\medskip

Thus for any large enough cardinal $\lambda$, there is a regular cardinal $\mu < \gimel_\omega$ such that $\clubsuit^\ast_{\lambda^+} [NS_{\lambda^+} \vert E^{\lambda^+}_\mu]$ holds. 


\medskip

Notice that by Facts 1.1 and 4.8 (iii), Fact 4.7 follows from Corollary 4.10 (ii).

\medskip

\begin{Obs} Let $\chi$ be an uncountable limit cardinal, and $\lambda$ a cardinal greater than or equal to $\chi$. Suppose that $\chi$ is a $\lambda$-revision cardinal. Then the following hold :
 \begin{enumerate}[\rm (i)] 
\item Assuming that $\alpha (\chi) \leq \lambda$, there is a cardinal $\sigma < \chi$ such that $\clubsuit^{cof \slash {\mu^+}, \ast}_{\lambda^+} [NS_{\lambda^+} \vert E^{\lambda^+}_\mu]$ holds for any regular cardinal $\mu$ with $\sigma \leq \mu < \chi$.
\item Assuming that $2^{< \chi} \leq \lambda$, there is a cardinal $\sigma < \chi$ such that $\clubsuit^\ast_{\lambda^+} [NS_{\lambda^+} \vert E^{\lambda^+}_\mu]$ holds for any regular cardinal $\mu$ with $\sigma \leq \mu < \chi$.
\item Assuming that $2^{< \chi} \leq \lambda$ and $2^\lambda = \lambda^+$, there is a cardinal $\sigma < \chi$ such that $\diamondsuit^\ast_{\lambda^+} [NS_{\lambda^+} \vert E^{\lambda^+}_\mu]$ holds for any regular cardinal $\mu$ with $\sigma \leq \mu < \chi$.
\end{enumerate}
\end{Obs}

{\bf Proof.} By Corollary 4.6 and Facts 2.2 and 4.8.
\hfill$\square$

\medskip


\medskip

\begin{Cor}  Let $\chi$ be an uncountable limit cardinal that is not a fixed point of the aleph function, and $\lambda$ a cardinal with $\chi < \lambda < \chi^{+ \chi}$. Then the following hold :
 \begin{enumerate}[\rm (i)] 
\item There is a cardinal $\sigma < \chi$ such that $\clubsuit^{cof \slash {\mu^+}, \ast}_{\lambda^+} [NS_{\lambda^+} \vert E^{\lambda^+}_\mu]$ holds for any regular cardinal $\mu$ with $\sigma \leq \mu < \chi$.
\item Assuming that $2^{< \chi} \leq \lambda$, there is a cardinal $\sigma < \chi$ such that $\clubsuit^\ast_{\lambda^+} [NS_{\lambda^+} \vert E^{\lambda^+}_\mu]$ holds for any regular cardinal $\mu$ with $\sigma \leq \mu < \chi$.
\item Assuming that $2^{< \chi} \leq \lambda$ and $2^\lambda = \lambda^+$, there is a cardinal $\sigma < \chi$ such that $\diamondsuit^\ast_{\lambda^+} [NS_{\lambda^+} \vert E^{\lambda^+}_\mu]$ holds for any regular cardinal $\mu$ with $\sigma \leq \mu < \chi$.
\end{enumerate}
\end{Cor}

{\bf Proof.} Use Proposition 2.7.
\hfill$\square$

\medskip

Notice that the main conclusion here is (i), since (ii) and (iii) will not apply in case $2^{< \chi} \geq \chi^{+ \chi}$.

\medskip

In Fact 4.7, we have, as in the RGCH, boundedly many \say{exceptions} below $\chi$. Finitely many would be a huge improvement, and this is what is attempted in \cite{SheMore}. In the words of Shelah in \cite{She10}, the result asserts that \say{for every $\lambda = \chi^+ = 2^\chi > \mu$, $\mu$ a strong limit for some \underline{finite} $
{\mathfrak d} \subseteq {\rm Reg} \cap \mu$, for every regular $\kappa < \mu$ not from ${\mathfrak d}$ we have $\diamondsuit_{S^\lambda_\kappa}$, and even $\diamondsuit_S$ for 
\say{most} stationary $S \subseteq S^\lambda_\kappa$. In fact, for the relevant good stationary sets $S \subseteq S^\lambda_\kappa$ we get $\diamondsuit_S^\ast$}. It remains to be 
seen whether a similar result holds for $\clubsuit^\ast$.

\section{An equality}

\bigskip

We will show that the RGCH has something to say concerning the two-place function $u (\sigma, \lambda)$.
Our study hinges on the following result of Shelah.

\medskip

\begin{fact} {\rm (\cite{SheFurther})} Let $\sigma, \chi, \lambda$ be three infinite cardinals with $\sigma \leq \chi < \lambda$. Then

\centerline{$u (\sigma, \nu) = \max\{u (\sigma, \chi), u (\chi^+, \lambda)\}$,}

where $\nu = \cov (\lambda, \chi^+, \chi^+, \sigma)$.
\end{fact}


\begin{Cor}  Let $\sigma < \chi < \lambda$ be three infinite cardinals. Suppose that $\cov (\lambda, \chi^+, \chi^+, \sigma) \leq \lambda$. Then

\centerline{$u (\sigma, \lambda) = \max\{u (\sigma, \chi), u (\chi^+, \lambda)\}$.}

\end{Cor}

\begin{Cor}  Let $\chi$ be an uncountable limit cardinal, and  $\lambda$ be a cardinal greater than $\chi$. Suppose that $\chi^+$ is a $\lambda$-revision cardinal. Then for any large enough cardinal $\sigma < \chi$, 

\centerline{$u (\sigma, \lambda) = \max\{u (\sigma, \chi), u (\chi^+, \lambda)\}$.}

\end{Cor}

{\bf Proof.} Use Fact 2.3 (ii).






\medskip

Under SSH, the picture is clear.

\medskip

\begin{fact} {\rm (\cite{Secret})} Let $\rho_1, \rho_2, \rho_3$ and $\rho_4$ be four infinite cardinals such that $\rho_1 \geq \rho_2 \geq \rho_3 \geq \rho_4$. Then assuming SSH, the following hold :        
\begin{enumerate}[\rm (i)]
\item  If $\rho_1 = \rho_2$ and either $\cf(\rho_1) < \rho_4$ or $\cf(\rho_1) \geq \rho_3$, then $\cov (\rho_1, \rho_2, \rho_3, \rho_4) = \cf(\rho_1)$.           
\item  If $\rho_4 \leq \cf (\rho_1) < \rho_3$, then $\cov (\rho_1, \rho_2, \rho_3, \rho_4) = \rho_1^+$.        
\item  In all other cases, $\cov (\rho_1, \rho_2, \rho_3, \rho_4) = \rho_1$.
 \end{enumerate}
\end{fact}

\begin{Cor}  Let $\sigma < \chi < \lambda$ be three infinite cardinals. Suppose that SSH holds and $\chi$ is a limit cardinal. Then the following are equivalent : 
\begin{enumerate}[\rm (i)]
\item Either $\cf (\lambda) < \sigma$, or $\chi < \cf (\lambda)$.
\item  $u (\sigma, \lambda) = \max\{u (\sigma, \chi), u (\chi^+, \lambda)\}$.
\end{enumerate}
\end{Cor}

\medskip




Let us look for more situations (in ZFC) when the conclusion of Corollary 5.3 holds. To make things easier we will only deal with the case when $\sigma$ is a regular cardinal. Even so, our results will require $\lambda$ to be close to $\chi$.


\medskip

\begin{fact} {\rm (\cite{DM})} Let $\kappa, \chi, \lambda$ be three uncountable cardinals with $\cf (\kappa) = \kappa \leq \chi < \lambda$. Then $u (\kappa, \lambda) \leq \max\{u (\kappa, \chi), u (\chi^+, \lambda)\}$.
\end{fact}

\medskip

For a cardinal $k$, $FP(k)$ denotes the least fixed point of the aleph function greater than $k$. 

\medskip

\begin{fact} {\rm (\cite{pcf})} Let $\kappa$ be an infinite successor cardinal. Then the following hold :
\begin{enumerate}[\rm (i)]
\item Let $\lambda$ be a cardinal with $\kappa < \lambda < \min\{FP(\kappa), u(\kappa, \lambda)\}$. Then $u (\kappa, \lambda)$ is less than $FP (\kappa)$ (and hence is not a weakly inaccessible cardinal). Moreover there is a (unique) cardinal $\theta (\kappa, \lambda)$ such that
 \begin{itemize}
 \item $\cf(\theta (\kappa, \lambda)) < \kappa < \theta (\kappa, \lambda) \leq \lambda$.
\item $u(\kappa, \lambda) = u (\kappa, \theta (\kappa, \lambda)) = \pp (\theta (\kappa, \lambda)) = \pp_{< \kappa} (\theta (\kappa, \lambda)) =$ 

$\cov(\theta (\kappa, \lambda), \theta (\kappa, \lambda), (\cf (\theta (\kappa, \lambda)))^+, \cf (\theta (\kappa, \lambda)))$. 
\item $\pp_{< \kappa} (\rho) < \theta (\kappa, \lambda)$ for any cardinal $\rho$ with $\cf(\rho) < \kappa \leq \rho < \theta (\kappa, \lambda)$.
 \item $\pp_{< \kappa} (\theta (\kappa, \lambda)) = \max \{ \pp_{< \kappa} (\rho) : \cf (\rho) < \kappa \leq \rho \leq \lambda\}$.
 \item If $\cf(\theta (\kappa, \lambda)) \not= \omega$, then $\pp(\theta (\kappa, \lambda)) = \pp_{\Gamma((\cf(\theta (\kappa, \lambda))^+, \cf(\theta (\kappa, \lambda)))} (\theta (\kappa, \lambda)) = \pp^{\ast}_{I_{\cf(\theta (\kappa, \lambda))}} (\theta (\kappa, \lambda))$.
\end{itemize}
\item $u (\kappa, FP (\kappa)) = \cov (FP (\kappa), FP (\kappa), \kappa, 2)$.
\end{enumerate}
\end{fact}

\begin{fact} {\rm (\cite[Theorem 5.4 p. 87]{SheCA})} Let $\nu, \theta$ be two infinite cardinals with $\cf (\theta) \leq \nu < \theta$. Then $\pp_\nu (\theta) \leq \cov(\theta, \theta, \nu^+,2)$.
\end{fact}

\begin{fact} {\rm (\cite{Fixed})} Let $\chi$ be a singular cardinal that is not a fixed point of the aleph function. Then for any large enough regular cardinal $\kappa < \chi$, $u (\kappa, \theta) < \chi$ for every cardinal $\theta$ with $\kappa \leq \theta < \chi$.
\end{fact}

\begin{Pro}  Let $\chi$ be an uncountable limit cardinal that is not a fixed point of the aleph function, and $\lambda$ be a cardinal with $\chi < \lambda < FP (\chi)$. Then letting $W$ denote the set of all regular cardinals less than $\chi$, the following hold :
\begin{enumerate}[\rm (i)]
\item One of the following holds :
\begin{enumerate}[\rm (a)]
\item For any large enough cardinal $\kappa$ in $W$, $u (\kappa, \lambda) = \lambda$.
\item There is a cardinal $\theta \geq \chi$ such that for any large enough cardinal $\kappa$ in $W$, ($u (\kappa, \lambda) > \lambda$ and) $\theta (\kappa, \lambda) = \theta$.
\end{enumerate}
\item For any large enough cardinal $\kappa$ in $W$, 

\centerline{$u (\kappa, \lambda) = \max\{u (\kappa, \chi), u (\chi^+, \lambda)\}$.}

\end{enumerate}
\end{Pro}

{\bf Proof.} 
Select a cardinal $\rho < \chi$ such that 
\begin{itemize}
\item $\lambda < FP (\rho)$.
\item $\cf (\theta (\chi^+, \lambda)) < \rho$ in case $u (\chi^+, \lambda) > \lambda$.
\end{itemize}
Then any cardinal $\tau$ with $\rho < \tau \leq \lambda$ is either a singular cardinal, or a successor cardinal.

(i) : 

\medskip  

{\bf Claim 1.} Let $\kappa, \mu \in W$. Suppose that $u (\kappa, \lambda) > \lambda$ and $\rho < \kappa < \mu < \theta (\kappa, \lambda)$. Then $u (\mu, \lambda) > \lambda$, and moreover $\theta (\mu, \lambda) \leq \theta (\kappa, \lambda)$.

\medskip

{\bf Proof of Claim 1.} By Fact 5.7, $u (\mu, \lambda) \geq u (\mu, \theta (\kappa, \lambda)) \geq \cov (\theta (\kappa, \lambda), \theta (\kappa, \lambda), \mu, 2)  \geq \cov (\theta (\kappa, \lambda), \theta (\kappa, \lambda), (\cf (\theta (\kappa, \lambda)))^+, 2) \geq \cov (\theta (\kappa, \lambda), \theta (\kappa, \lambda), (\cf (\theta (\kappa, \lambda)))^+, \cf (\theta (\kappa, \lambda))) = u (\kappa, \lambda) > \lambda$, so $u (\mu, \lambda) > \lambda$. Furthermore $\theta (\mu, \lambda) \leq \theta (\kappa, \lambda)$, since otherwise

\centerline{$u (\kappa, \lambda) = \pp (\theta (\kappa, \lambda)) \leq \pp_{< \mu} (\theta (\kappa, \lambda)) < \theta (\mu, \lambda) \leq \lambda$.}

This completes the proof of the claim.

\medskip

{\bf Claim 2.} One of the following holds :
\begin{enumerate}[\rm (1)]
\item (b).
\item For any large enough cardinal $\kappa$ in $W$, either $u (\kappa, \lambda) = \lambda$, or $u (\kappa, \lambda) > \lambda$ and $\theta (\kappa, \lambda) < \chi$.\end{enumerate}

\medskip

{\bf Proof of Claim 2.} Suppose that (2) does not hold. Then the set $X$ of all $\kappa \in W \setminus \rho^+$ such that $u (\kappa, \lambda) > \lambda$ and $\theta (\kappa, \lambda) \geq \chi$ is unbounded in $\chi$. Set $\theta = \min \{ \theta (\kappa, \lambda) : \kappa \in X \}$, and select $\tau \in X$ with $\theta (\tau, \lambda) = \theta$. Now let $\pi \in W \setminus \tau$. Then by Claim 1, $u (\pi, \lambda) > \lambda$. Furthermore $\theta (\pi, \lambda) \geq \chi$, since otherwise by Claim 1, for any $\eta \in W$ greater than $\pi$, ($u (\eta, \lambda) > \lambda$ and) $\theta (\eta, \lambda) \leq \theta (\pi, \lambda) < \chi$, contradicting the unboundedness of $X$. Thus $\pi \in X$. By a last appeal to Claim 1, $\theta \leq \theta (\pi, \lambda) \leq \theta (\tau, \lambda) = \theta$, which completes the proof of the claim.

\medskip

{\bf Claim 3.} Suppose that (2) of Claim 2 holds. Then (a) holds.

\medskip

{\bf Proof of Claim 3.} Suppose that (a) does not hold. Then the set of all $\kappa \in W \setminus \rho^+$ such that $u (\kappa, \lambda) > \lambda$ and $\theta (\kappa, \lambda) < \chi$ must be unbounded in $\chi$. This contradicts Fact 5.8,
which completes the proof of the claim and that of (i).

\medskip

(ii) : Let $\kappa$ be a regular cardinal with $\rho < \kappa < \chi$. By Fact 5.5, 

\centerline{$u (\kappa, \lambda) \leq \max\{u (\kappa, \chi), u (\chi^+, \lambda)\}$.}

In fact, we have equality. Clearly, $u (\kappa, \lambda) \geq u (\kappa, \chi)$, so we need to show that $u (\kappa, \lambda) \geq u (\chi^+, \lambda)$. This is immediate in case $u (\chi^+, \lambda) = \lambda$. Otherwise by Fact 5.7, letting $\kappa = \nu^+$, we have that $u (\kappa, \lambda) = u (\nu^+, \lambda) \geq u (\nu^+, \theta (\chi^+, \lambda)) \geq \cov (\theta (\chi^+, \lambda), \theta (\chi^+, \lambda), \nu^+, 2)  \geq \pp_\nu (\theta (\chi^+, \lambda)) \geq  \pp (\theta (\chi^+, \lambda)) = u (\chi^+, \lambda)$.

\medskip

\hfill$\square$

\begin{Pro}  Let $\chi$ be a fixed point of the aleph function that is not a limit of fixed points of the aleph function, and $\lambda$ be a cardinal with $\chi < \lambda < FP (\chi)$. Then for any large enough regular cardinal $\kappa < \chi$,

\centerline{$u (\kappa, \lambda) = \max\{u (\kappa, \chi), u (\chi^+, \lambda)\}$.}
\end{Pro}

{\bf Proof.} The proof is a slight modification of that of Proposition 5.10 (ii). 
Notice that $\cf (\chi) = \omega$.
Select a cardinal $\rho < \chi$ such that 
\begin{itemize}
\item $FP (\rho) = \chi$.
\item $\cf (\theta (\chi^+, \lambda)) < \rho$ in case $u (\chi^+, \lambda) > \lambda$.
\end{itemize}
Then any regular cardinal $\tau$ with $\rho < \tau \leq \chi$ is a successor cardinal. Now let $\kappa$ be a regular cardinal with $\rho < \kappa < \chi$. 

\medskip

$u (\kappa, \lambda) \leq \max\{u (\kappa, \chi), u (\chi^+, \lambda)\}$ : By Fact 5.5. 

\medskip

$u (\kappa, \lambda) \geq u (\kappa, \chi)$ : Clear.

\medskip

$u (\kappa, \lambda) \geq u (\chi^+, \lambda)$ : This is immediate in case $u (\chi^+, \lambda) = \lambda$. Otherwise by Fact 5.7, letting $\kappa = \nu^+$, we have that $u (\kappa, \lambda) = u (\nu^+, \lambda) \geq u (\nu^+, \theta (\chi^+, \lambda)) \geq \cov (\theta (\chi^+, \lambda), \theta (\chi^+, \lambda), \nu^+, 2)  \geq \pp_\nu (\theta (\chi^+, \lambda)) \geq  \pp (\theta (\chi^+, \lambda)) = u (\chi^+, \lambda)$.
\hfill$\square$

\bigskip

\section{The virtue of dishonesty}

\bigskip


The birth of the RGCH is described in \cite{SheCard} as follows : \say{Condition (b) of Theorem 2 holds easily for $\mu = \lambda$. Still it may look restrictive, and the author was tempted to try to eliminate it (...). But instead of working \say{honestly} on this the author for this purpose proved (see [She460]) that it follows from ZFC, and therefore can be omitted (...)}, where [She460] = \cite{SheRGCH}. Let us note that the Theorem 2 in question is a topological statement, so the RGCH is being applied from the very beginning.
\bigskip

\bigskip

{\bf Acknowledgements}. The author would like to thank Assaf Rinot for pointing out a mistake in a previous version of the article.

  \bigskip
\noindent Universit\'e de Caen - CNRS \\
Laboratoire de Math\'ematiques \\
BP 5186 \\
14032 Caen Cedex\\
France\\
Email :  pierre.matet@unicaen.fr\\


\end{document}